\documentclass[10pt]{article}
\usepackage{amsfonts}
\usepackage{amsmath}
\usepackage{mathrsfs}
\usepackage{mathrsfs,amscd,amssymb,amsthm,amsmath,bm,graphicx,psfrag,subfigure}

\makeatletter

\renewcommand{\@seccntformat}[1]{{\csname the#1\endcsname}{\normalsize .}\hspace{.5em}}
\makeatother

\def \[{\begin{equation}}
\def \]{\end{equation}}

\newtheorem{thm}{Theorem}[section]

\newtheorem{claim}{Claim}
\newtheorem{fact}{Fact}
\newtheorem{lem}[thm]{Lemma}

\newenvironment{wst}
{\setlength{\leftmargini}{1.5\parindent}
 \begin{itemize}
 \setlength{\itemsep}{-1.1mm}}
{\end{itemize}}

\begin{document}
\setlength{\baselineskip}{15pt}
\begin{center}{\Large \bf CORRIGENDUM\\[3pt] `On the signless Laplacian spectra of $k$-trees'-- CORRIGENDUM\footnote{Financially supported by the National Natural
Science Foundation of China (Grant Nos. 11271149,11371162), the
Program for New Century Excellent Talents in University (Grant No.
NCET-13-0817) and the Special Fund for Basic Scientific Research of
Central Colleges (Grant No. CCNU13F020).}}

\vspace{5mm}

{\large Minjie Zhang$^{a,b,}$\footnote{Corresponding author. \\
\hspace*{5mm} E-mail:
zmj1982@21cn.com(M.J. Zhang), lscmath@mail.ccnu.edu.cn (S.C. Li)}, Shuchao Li$^{a}$}\vspace{2mm}

$^a$Faculty of Mathematics and Statistics,  Central China Normal
University, Wuhan 430079, P.R. China

$^b$Faculty of Mathematics and Physics, Hubei Institute of
Technology, Huangshi 435003, P.R. China
\end{center}
\vspace{4mm}

A $k$-\textit{tree} is either a complete graph on $k$ vertices or a graph obtained from a smaller $k$-tree by adjoining a new vertex together with $k$ edges connecting it to a $k$-clique. We tried to determine the $n$-vertex $k$-trees with the first three largest signless Laplacian indices in [Linear Algebra Appl.  467 (2015) 136-148]. Our proofs were mainly based on Lemma 2.4 in the above paper, whereas Lemma 2.4 follows from the next lemma:

\noindent($\bullet$) (Lemma 2.3 of \cite{Z-B-X})\textit{
Let $0<a_i\leqslant 1, 0<b_i\leqslant 1$ for $1\leqslant  i\leqslant  k$. Then
$
\sum_{i=1}^kb_i^2-\sum_{i=1}^ka_i^2 \leqslant 2\left(\sum_{i=1}^kb_i-\sum_{i=1}^ka_i\right)
$
and the equality holds if and only if $a_i=b_i=1$ for $i=1,2,\ldots,k$.}

In fact, in the proof of Lemma 2.3, if $b_{i}-a_{i}<0$,  one cannot conclude the result. For example if $\sum_{1}^{k}b_{i}=\sum_{1}^{k}a_{i}$, we cannot conclude $\sum_{1}^{k}b_{i}^{2}-\sum_{1}^{k}a_{i}^{2}\leqslant 0$. Hence, Lemma 2.3 is not correct, which implies that our method in \cite{Z-B-X} can not be used to determine the $n$-vertex $k$-trees with the first three largest signless Laplacian indices.

The text below give a new method to determine the $n$-vertex $k$-trees with the first three largest signless Laplacian indices.

% {\setcounter{section}{0}
%\section{\normalsize Introduction}\setcounter{equation}{0}

\section{\normalsize Preliminaries}
%In this section, we give some preliminary results.
%\begin{lem}[\cite{X-M-J}]
%Let $T_n^k$ be a $k$-tree on $n\geqslant  k+1$ vertices. Then, there exists at least two vertices in $V_{T_n^k}$ of degree $k$.
%\end{lem}
%\begin{lem}[\cite{M-B}]
%Let $G$ be a connected graph with at least one edge and maximum degree $\Delta$, then $q_1(G)\geqslant \Delta(G)+1$ with equality
%if and only if $G\cong S_n$, where $S_n$ is the star graph on $n$ vertices.
%\end{lem}
%In this section, we give some necessary definitions and some basic facts on $n$-vertex $k$-trees.
We only consider simple connected graphs $G=(V_G, E_G)$, where $V_G$ is the vertex set and $E_G$ is the edge set. We call $n=|V_G|$ the \textit{order} of $G$ and $m=|E_G|$ the \textit{size} of $G$. For a vertex subset $S$ of $V_G$, denoted by $G[S]$ the subgraph induced by $S$. Let $N_G(v)$ denote the set of vertices adjacent to $v.$ $N_G[v]=N_G(v)\cup\{v\}.$ Let $A(G)$ be the adjacency matrix of graph $G$, $D(G)$ be the diagonal matrix with degrees of the vertices on the main diagonal. The signless matrix $Q(G)=D(G)+A(G)$ is real symmetric, its eigenvalues can be arranged as $q_1(G)\geqslant q_2(G)\geqslant \cdots \geqslant q_n(G)\geqslant 0$, where $q_1(G)$ is the \textit{signless Laplacian index} of graph $G$.
\begin{figure}[h!]
\begin{center}
  % Requires \usepackage{graphicx}
\psfrag{1}{$K_k$}\psfrag{2}{$G_1$}
\psfrag{3}{$v_1$}\psfrag{4}{$v_{k-1}$}\psfrag{5}{$v_k$}
\psfrag{t}{$S_{k,n-k}$}\psfrag{x}{$u_3$}\psfrag{w}{$u_4$}\psfrag{r}{$u_{n-k-1}$}
\psfrag{u}{$u_1$}\psfrag{v}{$u_2$}\psfrag{y}{$u_{n-k}$}
\includegraphics[width=70mm]{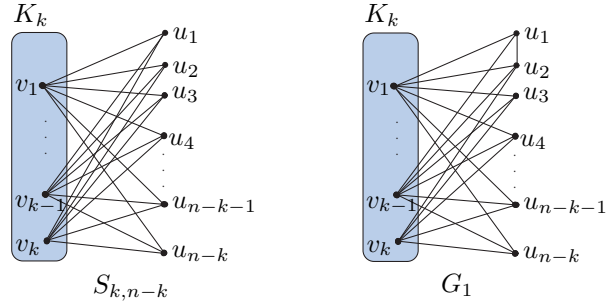}\\
  \caption{Graphs $S_{k,n-k}$ and $G_1$.}\label{figure 1}
\end{center}
\end{figure}

For positive integers $n,k$ with $n\geqslant  k$, the $k$-\textit{star} $S_{k,n-k}$ is depicted in Fig. 1. Let $G_1$ be the graph obtained from $S_{k,n-k}$ by deleting the edge $u_1v_k$ and adding an edge $u_1u_2$ (see Fig. 1). For convenience, let $\mathscr{T}_n^k$ denote the set of all $k$-trees on $n$ vertices. Let $T_n^k$ be an $n$-vertex $k$-tree. Obviously, $T_n^k\cong K_k$ if $n=k$. So we consider $T_n^k$ with $n\geqslant  k+1$ in the whole context.  A vertex $v\in V_{T_n^k}$ is called a $k$-\textit{simplicial vertex} if $v$ is a vertex of degree $k$ whose neighbors form a $k$-clique of $T_n^k$.
Let $S_1(T_n^k)$ be the set of all the $k$-simplicial vertices of $T_n^k$ and we can easily have the following facts:
\begin{fact}
$S_1(T_n^k)=\emptyset$ if and only if $T_n^k\cong K_k$; For any $v$ in $V_{K_{k+1}}$,
$v$ is in $S_1(K_{k+1})$; If $n\geqslant    k+2$, then $|S_1(T_n^k)|\geqslant    2$ and $S_1(T_n^k)$ is
an independent set.
\end{fact}
\begin{fact}
Let $T_n^k$ be a $k$-tree on $n\geqslant k+2$ vertices. Then $d_{T_n^k}(v)=k$ for any $v\in S_1(T_n^k)$ and $T_n^k-v$ is
also a $k$-tree.
\end{fact}
\begin{fact}
$|S_1(T_n^k)|=n-k$ if and only if $T_n^k\cong S_{k,n-k}$.
\end{fact}

Given a $k$-clique $G[\{v_1,v_2,\ldots,v_k\}]$, if there exists another vertex $w$ in $G$ such that $G[\{v_1,v_2,\ldots,v_k,w\}]$ is a $(k+1)$-clique and $d_G(w)=k$, then we say that $G[\{v_1,v_2,\ldots,v_k,w\}]$ has the \textit{property} $P_G(v_1,v_2,\ldots,v_k)$.  Let $l_G(v_1,v_2,\ldots,v_k)$ be the number of $(k+1)$-cliques which has the property $P_G(v_1,v_2,\ldots,v_k)$ and set
$$
  \text{$l(G):=\max_{\{v_1,\ldots, v_k\}\subseteq V_G}\{l_G(v_1,v_2,\ldots,v_k): G[\{v_1, v_2,\ldots, v_k\}]$ is a $k$-clique of $G$\}}.
$$

Obviously, $l(G)=n-k$ if and only if $G\cong S_{k,n-k}$, whereas $l(G)=n-k-2$ if and only if $G\cong G_1$; and there does not exist an $n$-vertex $k$-tree, say $G$, with $l(G)=n-k-1$.
\begin{fact}
Let $G$ be an $n$-vertex $k$-tree. Then $l(G)=n-k-3$ if and only if $G\cong G_2, G_3, G_4$ or $G_5$, where $G_2=G_1-u_3v_k+u_3u_2, G_3=G_1-u_3v_k+u_3u_1,\, G_4=G_1-u_3v_1+u_3u_2,\, G_5=G_1-u_3v_{k-1}-u_3v_k+u_3u_1+u_3u_2$. Graphs $G_2, G_3, G_4, G_5$ are depicted in Fig. 2.
\end{fact}
\begin{fact}
For graphs $G_2, G_3, G_4, G_5$ (see Fig. 2), one has $|S_1(G_2)|=|S_1(G_4)|=n-k-1,\, |S_1(G_3)|=|S_1(G_5)|=n-k-2$.
\end{fact}
\begin{figure}[h!]
\begin{center}
 % Requires \usepackage{graphicx}
\psfrag{1}{$K_k$}\psfrag{2}{$n-k-3$}
\psfrag{3}{$v_1$}\psfrag{4}{$v_{k-1}$}\psfrag{5}{$v_k$}
\psfrag{u}{$u_1$}\psfrag{v}{$u_2$}\psfrag{x}{$u_3$}\psfrag{y}{$u_{n-k}$}
\psfrag{6}{$G_2$}\psfrag{7}{$G_3$}\psfrag{w}{$u_4$}\psfrag{t}{$u_{n-k-1}$}
\psfrag{8}{$G_4$}\psfrag{9}{$G_5$}
\includegraphics[width=130mm]{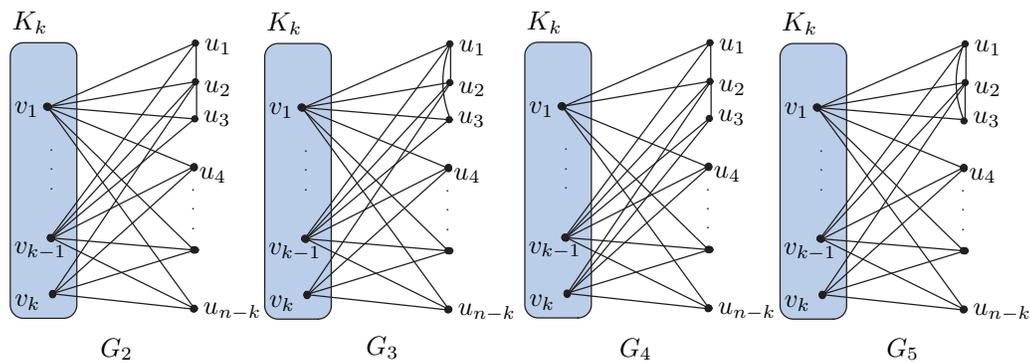}\\
  \caption{Graphs $G_2,G_3,G_4$ and $G_5$.}\label{figure 3}
\end{center}
\end{figure}

Further on we will need the following lemma.
\begin{lem}[\cite{H-Z}]
Let $u$ and $v$ be two distinct vertices of a connected graph $G$. Suppose that $w_1,w_2,\ldots,w_{s}\, (s\geqslant 1)$ are neighbors of
$v$ but not $u$ and they are all different from $u$. Let ${\bf x}=(x_1, x_2, \ldots, x_n)^T$ be the Perron vector of $Q(G)$, and
let $G^*$ be obtained from $G$ by deleting the edges $vw_{i}$ and adding the edges $uw_{i}$ for $i = 1,2,\ldots,s$. If $x_v \leqslant
x_u,$ then $q_1(G) < q_1(G^*)$.
\end{lem}

\section{\normalsize Characterization of the $k$-trees with the first three largest signless Laplacian indices}\setcounter{equation}{0}
In this section, we characterize the $k$-trees on $n$ vertices having the maximal, second maximal and third maximal signless Laplacian indices.
\begin{lem}
Let $G$ be in $\mathscr{T}_n^k\setminus \{S_{k,n-k}\}$ on $n\geqslant  k+1$ vertices. Then there exists a $k$-tree $G^*$ on $n$ vertices
such that $q_1(G) < q_1(G^*)$ with $|S_1(G^*)|=|S_1(G)|+1$.
\end{lem}
\begin{proof}
Note that $G\ncong S_{k,n-k}$. By Fact 3, we get $|S_1(G)|<n-k$. Then, $G-S_1(G)$ is a $k$-tree whose order is no less than $k+1$ and $S_1(G-S_1(G))\neq \emptyset$ by Facts 1 and 2. Let $u$ be a vertex in $S_1(G-S_1(G))$ with $N_{G-S_1(G)}(u)=\{v_1,v_2\ldots, v_k\}$. Assume that $N_G(u)\cap S_1(G)=\{w_1,w_2,\ldots,w_s\}$. Then, $d_G(w_j)=k$ and $N_G(w_j)\subseteq \{v_1,v_2,\ldots,v_k,u\}$ for any $1\leqslant j\leqslant s$ by the definition of the $k$-tree. In order to complete the proof, it suffices to show the following two claims.
\begin{claim}
If $G$ contains two vertices $w_i$ and $w_j$ in $\{w_1,w_2,\ldots,w_s\}$ satisfying $N_G(w_i)\neq N_G(w_j)$, then there exists an $n$-vertex $k$-tree $G'$ such that $q_1(G) < q_1(G')$ and $N_{G'}(w_i)= N_{G'}(w_j)$ with $|S_1(G')|=|S_1(G)|$ or $|S_1(G)|+1$.
\end{claim}
\noindent{\bf Proof of Claim 1.}\
Without loss of generality, we assume
$
  N_G(w_i)=\{v_1,\ldots,v_k,u\}\setminus\{v_i\},\ N_G(w_j)=\{v_1,\ldots,\linebreak v_k,u\}\setminus\{v_j\}
$
and
$$
  l_G(v_1,\ldots,v_{i-1},v_{i+1},\ldots,v_k,u)=a,\ \ \ \   l_G(v_1,\ldots,v_{j-1},v_{j+1},\ldots,v_k,u)=b,
$$
where $v_i\neq v_j$. Hence, $1\leqslant    a\leqslant    s,\ 1\leqslant    b\leqslant    s$. Let $\delta_1,\delta_2,\ldots,\delta_a$ (resp. $\Delta_1,\Delta_2,\ldots,\Delta_b$) be the
$(k+1)$-cliques with the property $P_G(v_1,\ldots,v_{i-1},v_{i+1},\ldots,v_k,u)$ (resp. $P_G(v_1,\ldots,v_{j-1},v_{j+1},\ldots,v_k,u)$) and
$w_{i_1}(=w_i),w_{i_2},\ldots,w_{i_a}$ (resp. $w_{j_1}(=w_j),w_{j_2},\ldots,w_{j_b}$) be the vertices in these $(k+1)$-cliques
different from $v_1,\ldots,v_{i-1},v_{i+1},\ldots,v_k,u$ (resp. $v_1,\ldots,v_{j-1},v_{j+1},\ldots,v_k,u$). Suppose, without loss of generality,
that $x_{v_i}\geqslant x_{v_j}$. Let
$$
  G'=G-\{v_jw_{i_1},v_jw_{i_2},\ldots,v_jw_{i_a}\}+\{v_iw_{i_1},v_iw_{i_2},\ldots,v_iw_{i_a}\}.
$$
Then we have $G'\in \mathscr{T}_n^k$ with $N_{G'}(w_i)= N_{G'}(w_j)$ and $q_1(G) < q_1(G')$ by Lemma 1.1. Note that
$$
  \{v_1,\ldots,v_{j-1},v_{j+1},\ldots,v_k,u,w_{i_1},w_{i_2},\ldots,w_{i_a}\}\subseteq N_G(v_j).
$$
Hence, we get $S_1(G')=S_1(G)$ if $d_G(v_j)> k+a$; whereas $S_1(G')=S_1(G)\cup \{v_j\}$ if $d_G(v_j)= k+a$. Thus, $|S_1(G')|=|S_1(G)|$ or $|S_1(G)|+1$, as desired.
\qed

\begin{claim}
If $N_G(w_1)= N_G(w_2)=\cdots = N_G(w_s)$, then there exists an $n$-vertex $k$-tree $G^*$ such that $q_1(G) < q_1(G^*)$ and $|S_1(G^*)|=|S_1(G)|+1$.
\end{claim}
\noindent{\bf Proof of Claim 2.}\
Without loss of generality, we assume $N_G(w_j)=\{v_1,\ldots,v_{k-1},u\}$ for $1\leqslant j\leqslant s$. Let $V_0=N_G(v_k)\setminus \{v_1,\ldots,v_{k-1},u\}$. Then, $V_0\neq\emptyset$ since $v_k\not\in S_1(G)$. If $x_{v_k}\geqslant x_u$, then let
$$
   G^*=G-\{uw_1,uw_2,\ldots,uw_s\}+\{v_kw_1,v_kw_2,\ldots,v_kw_s\}.
$$
It is routine to check that $G^*$ is an $n$-vertex $k$-tree and $|S_1(G^*)|=|S_1(G)\cup \{u\}|=|S_1(G)|+1.$ By Lemma 1.1, $q_1(G) < q_1(G^*)$.

If $x_{v_k}< x_u$, then let $G^*=G-\sum_{v\in V_0}v_kv+\sum_{v\in V_0}uv$. It is routine to check $G^*$ is an $n$-vertex $k$-tree and $|S_1(G^*)|=|S_1(G)\cup \{v_k\}|=|S_1(G)|+1$. By Lemma 1.1, $q_1(G) < q_1(G^*)$.

This completes the proof of Claim 2.
\qed

Now we come back to the proof of Lemma 2.1. If $N_G(w_1)= N_G(w_2)=\cdots = N_G(w_s)$, by Claim 2 our result holds. Otherwise, repeatedly using Claim 1 (at most $(s-1)$ times), one can finally get an $n$-vertex $k$-tree $G''$ with $q_1(G) < q_1(G'')$ satisfying one of the following:
\begin{wst}
\item[{\rm (i)}] $|S_1(G'')|=|S_1(G)|+1$;
\item[{\rm (ii)}] $|S_1(G'')|=|S_1(G)|$ and $N_{G''}(w_1)= N_{G''}(w_2)=\cdots = N_{G''}(w_s)$.
\end{wst}
If (i) holds, then let $G^*\cong G''$, our lemma holds obviously. If (ii) holds, combining with Claim 2, we can get a $k$-tree $G^*$ on $n$ vertices
such that $q_1(G'') < q_1(G^*)$ and $|S_1(G^*)|=|S_1(G'')|+1=|S_1(G)|+1$, as desired.
\end{proof}

The following theorem follows immediately from Lemma 2.1.
\begin{thm}
Let $G$ be an $n$-vertex $k$-tree with $n\geqslant k+1$. Then $q_1(G)\leqslant q_1(S_{k,n-k})$ and the equality holds if and only if $G\cong
S_{k,n-k}$, where $S_{k,n-k}$ is depicted in Fig. 1.
\end{thm}
\begin{lem}
Let $G$ be an $n$-vertex $k$-tree with $n\geqslant  k+1$ and $|S_1(G)|=n-k-1$. If $G\ncong G_1$, then there exists an $n$-vertex $k$-tree $G^*$ with $|S_1(G^*)|=n-k-1$ such that  $q_1(G) < q_1(G^*)$ and $l(G^*)=l(G)+1$.
\end{lem}
\begin{proof}
Note that $|S_1(G)|=n-k-1$. Hence, the order of $G-S_1(G)$ is $k+1$. By Fact 2, $G-S_1(G)$ is a $k$-tree, that is, $G-S_1(G)\cong K_{k+1}$. Denote the vertex set of $G-S_1(G)$ by $\{v_1,v_2,\ldots,v_{k+1}\}$.

As $G\ncong G_1$, we have $l(G)<n-k-2$. Assume, without loss of generality, that $l(G)=l_G(v_1,v_2,\ldots,v_k)=s$. Hence, $s< n-k-2$. Let $\delta_1,\delta_2,\ldots,\delta_s$ be the $(k+1)$-cliques satisfying the property $P_G(v_1,v_2,\ldots,v_k)$ and
$w_1,w_2,\ldots,w_s$ be the vertices in these $(k+1)$-cliques different from $v_1,v_2,\ldots,v_k$. Obviously,
$d_G(w_1)=d_G(w_2)=\cdots=d_G(w_s)=k$ and $w_iw_j\not\in E_G$ for $1\leqslant i< j\leqslant s$.

Note that $|S_1(G)|=n-k-1$ and $s< n-k-2$. Hence, there exists a vertex $r_1\not\in \{w_1,w_2,\ldots,w_s\}$ such that $d_G(r_1)=k$ and
$G\left[N_G[r_1]\right]\cong K_{k+1}$. Without loss of generality, we assume $N_G(r_1)=\{v_2,v_3,\ldots,v_{k+1}\}$ and $l_G(v_2,v_3,\ldots,v_{k+1})=t$.
Obviously, $1\leqslant t\leqslant s$. Let $\Delta_1,\Delta_2,\ldots,\Delta_t$ be the $(k+1)$-cliques satisfying the property $P_G(v_2,v_3,\ldots,v_{k+1})$ and
$r_1,r_2,\ldots,r_t$ be the vertices in these $(k+1)$-cliques different from $v_2,v_3,\ldots,v_{k+1}$. Suppose that ${\bf x}=(x_1, x_2,
\ldots, x_n)^T$ is the Perron vector of $Q(G)$.% We proceed by distinguishing the following two possible cases:\vspace{2mm}

If $x_{v_1}\geqslant x_{v_{k+1}}$, then let $G^*=G-v_{k+1}r_1+v_1r_1$. It is routine to check that $G^*\in \mathscr{T}_n^k$. By Lemma 1.1, $q_1(G) < q_1(G^*).$
As $s+1< n-k-1=|S_1(G)|$, there exists a vertex $y\in S_1(G)\setminus \{w_1,w_2,\ldots,w_s,r_1\}.$ Note that $G[N_G(y)]$ is a $k$-clique not satisfying the property $P_G(v_1,v_2,\ldots,v_k)$. Hence, $yv_{k+1}\in E_G$. Thus, $v_{k+1}\not\in S_1(G^*)$. We have $l(G^*)=|\{w_1,w_2,\ldots,w_s,r_1\}|=l(G)+1$ and $S_1(G^*)=S_1(G)$, as desired.

If $x_{v_1}< x_{v_{k+1}}$, then let
$$
  G^*=G-\{v_1w_1,v_1w_2,\ldots,v_1w_{s-t+1}\}+\{v_{k+1}w_1,v_{k+1}w_2,\ldots,v_{k+1}w_{s-t+1}\}.
$$
Then $G^*$ is in $\mathscr{T}_n^k$ and by Lemma 1.1, $q_1(G) < q_1(G^*)$. Note that if $t\geqslant 2$, then one may easily get $v_1\not\in S_1(G^*)$. Now we consider that $t=1$. In this case, as $s+1< n-k-1=|S_1(G)|$, there exists a vertex $z\in S_1(G)\setminus \{w_1,w_2,\ldots,w_s,r_1\}$. Note that $G[N_G(z)]$ is a $k$-clique not satisfying the property $P_G(v_2,v_3,\ldots,v_{k+1})$. Hence, $zv_1\in E_G$ and $v_1\not\in S_1(G^*)$. Thus,  $l(G^*)=|\{r_1,r_2,\ldots,r_t,w_1,w_2,\ldots,w_{s-t+1}\}|=s+1=l(G)+1$ and $S_1(G^*)=S_1(G)$, as desired.
\end{proof}

The following theorem follows immediately from Lemmas 2.1 and 2.3.
\begin{thm}
Let $G$ be in $\mathscr{T}_n^k\setminus \{S_{k,n-k}\}$ with $n\geqslant k+1$. Then $q_1(G)\leqslant  q_1(G_1)$ with equality if
and only if $G\cong G_1,$ where $G_1$ is depicted in Fig. 1.
\end{thm}
\begin{lem}
Let $G$ be an $n$-vertex $k$-tree with $n\geqslant  k+1$ and $l(G)=n-k-3$. Then $q_1(G)\leqslant  q_1(G_2)$ with equality if and only if $G\cong
G_2$, where $G_2$ is depicted in Fig. 2.
\end{lem}
\begin{proof}
As $G$ is an $n$-vertex $k$-tree with $l(G)=n-k-3$, by Fact 4, it suffices to show the next three claims.
\setcounter{claim}{0}
\begin{claim}
$q_1(G_3) < q_1(G_2)$.
\end{claim}
\noindent{\bf Proof of Claim 1.}\
Let ${\bf x}=(x_1, x_2, \ldots, x_n)^T$ be the Perron vector of $Q(G_3)$. Put $G':=G_3-u_1u_3+u_2u_3$ if $x_{u_1}\leqslant x_{u_2}$ and
$G'':=G_3-u_2v_k+u_1v_k$ otherwise. It is routine to check that $G'\cong G''\cong G_2$. By Lemma 1.1, we have $q_1(G_3) < q_1(G_2),$ as desired.
\qed

By the same discussion as in the proof of Claim 1 as above, we can show the next claim, which is omitted.% here.
\begin{claim}
$q_1(G_4) < q_1(G_2)$.
\end{claim}
%\noindent{\bf Proof of Claim 2.}\
%Let ${\bf x}=(x_1, x_2, \ldots, x_n)^T$ be the Perron vector of $Q(G_4)$. If $x_{v_1}\geqslant x_{v_k}$, let $G'=G_4-v_ku_3+v_1u_3$;
%otherwise, let $G''=G_4-v_1u_1+v_ku_1$. It is straightforward to check that $G'\cong G''\cong G_2.$  By Lemma 1.1, we have $q_1(G_4) < q_1(G_2),$ as desired.
%\qed

\begin{claim}
$q_1(G_5) < q_1(G_2)$.
\end{claim}
\noindent{\bf Proof of Claim 3.}\ Let ${\bf x}=(x_1, x_2, \ldots, x_n)^T$ be the Perron vector of $Q(G_5)$. If $x_{v_{k-1}}\geqslant x_{u_2}$, then let $G':=G_5-u_2u_3+v_{k-1}u_3.$ Clearly, $G'\cong G_3$. By Lemma 1.1 and Claim 1, one has $q_1(G_5) < q_1(G_3) < q_1(G_2)$, as desired.

Now we consider $x_{v_{k-1}}<x_{u_2}$. For convenience, let $V_0=N_{G_5}(v_{k-1})\setminus \{v_1,v_2,\ldots,v_{k-2},v_k,u_1,u_2\}$. In this case, put $G'':=G_5-\sum_{v\in V_0}vv_{k-1}+\sum_{v\in V_0}vu_2$. It is routine to check that $G''\cong G_3$.  By Lemma 1.1 and Claim 1, one has $q_1(G_5) < q_1(G_3) < q_1(G_2)$, as desired.
\qed

By Claims 1-3, Lemma 2.5 holds.
\end{proof}

\begin{lem}
Let $G$ be an $n$-vertex $k$ tree with $n\geqslant k+1$ and $|S_1(G)|=n-k-2$. If $l(G)<n-k-3$, then there exists an $n$-vertex $k$ tree $G^*$ with $l(G^*)=l(G)+1$ such that $q_1(G) < q_1(G^*)$ and one of the following holds:
\begin{eqnarray}
|S_1(G^*)|&=&n-k-2,\\
|S_1(G^*)|&=& n-k-1\ and\ G^*\cong G_2.
\end{eqnarray}
\end{lem}
\begin{proof}
Note that $|S_1(G)|=n-k-2$; hence, together with Fact 2, $G-S_1(G)$ is a $(k+2)$-vertex $k$-tree, i.e., $G-S_1(G)\cong S_{k,2}$. Let $V_{G-S_1(G)}=\{v_1,v_2,\ldots,v_{k+2}\}$ with $d_{G-S_1(G)}(v_{k+1})=d_{G-S_1(G)}(v_{k+2})=k$. Obviously, $v_{k+1}v_{k+2}\not\in E_G$. Clearly, $l(G)\geqslant l_G(v_1,v_2,\ldots,v_k)$. Put $l_G(v_1,v_2,\ldots,v_k)=s$. %\vspace{2mm}

{\bf Case 1.}\ $l(G)=l_G(v_1,v_2,\ldots,v_k)$. In this case, one has $s< n-k-3$. Let $\delta_1,\delta_2,\ldots,\delta_s$ be the $(k+1)$-cliques satisfying the property $P_G(v_1,v_2,\ldots,v_k)$ and $w_1,w_2,\ldots,w_s$ be the vertices in these $(k+1)$-cliques different from $v_1,v_2,\ldots,v_k$. Then by Fact 1,
$d_G(w_1)=d_G(w_2)=\cdots=d_G(w_s)=k$ and $w_iw_j\not\in E_G$ for $1\leqslant i< j\leqslant s.$

Obviously, there exists a vertex $r_1$ in $S_1(G)$ such that $v_{k+1}r_1\in E_G$ by the definitions of the $k$-tree and $S_1(G)$. Assume, without loss of generality, that $N_G(r_1)=\{v_2,v_3,\ldots,v_k,v_{k+1}\}$ and $l_G(v_2,v_3,\ldots,v_k,v_{k+1})=t$. Then, $1\leqslant t\leqslant s$. Let $\Delta_1,\Delta_2,\ldots,\Delta_t$ be the $(k+1)$-cliques satisfying the property $P_G(v_2,v_3,\ldots,v_k,v_{k+1})$ and
$r_1,r_2,\ldots,r_t$ be the vertices in these $(k+1)$-cliques different from $v_2,v_3,\ldots,v_k,v_{k+1}$. Suppose that ${\bf x}=(x_1, x_2,
\ldots, x_n)^T$ is the Perron vector of $Q(G)$.

Since $v_{k+1}, v_{k+2}\not\in S_1(G)$, we get $d_G(v_{k+1}),\, d_G(v_{k+2})\geqslant k+1.$ We proceed by considering the following two possible subcases:\vspace{2mm}

\textbf{Subcase 1.1.}\ There exists a vertex, say $v_{k+1}$, in $\{v_{k+1},v_{k+2}\}$ such that $d_G(v_{k+1})> k+1$. If $x_{v_1}\geqslant x_{v_{k+1}}$, then let $G^*=G-v_{k+1}r_1+v_1r_1$. By Lemma 1.1, we have $q_1(G) < q_1(G^*).$  Combining with $d_G(v_{k+1})> k+1$, one can easily check that $l(G^*)=l(G)+1$ and $S_1(G^*)=S_1(G)$. If  $x_{v_1} < x_{v_{k+1}}$, then let
$$
  G^*=G-\{v_1w_1,v_1w_2,\ldots,v_1w_{s-t+1}\}+\{v_{k+1}w_1,v_{k+1}w_2,\ldots,v_{k+1}w_{s-t+1}\}.
$$
Then $G^*$ is in $\mathscr{T}_n^k$ and by Lemma 1.1, $q_1(G) < q_1(G^*)$. Note that $\{v_2,v_3,\ldots,v_k,v_{k+1},v_{k+2}\}\subseteq N_{G^*}(v_1)$,  we also get $l(G^*)=l(G)+1$ and $S_1(G^*)=S_1(G)$, our Lemma holds in this case.

\textbf{Subcase 1.2.}\ $d_G(v_{k+1})=d_G(v_{k+2})= k+1$. In this case, we have $s=l(G)=n-k-4$. Without loss of generality, we assume $x_{v_{k+1}}\leqslant x_{v_{k+2}}$ and let $G'=G-v_{k+1}r_1+v_{k+2}r_1$. By Lemma 1.1, we have $q_1(G) < q_1(G')$. Combining with $d_G(v_{k+1})= k+1$, we get
$$
  l(G')=|\{w_1,w_2,\ldots,w_s,v_{k+1}\}|=l(G)+1=n-k-3.
$$
However, noting that $S_1(G')=S_1(G)\cup\{v_{k+1}\}$, hence, $|S_1(G')|=|S_1(G)|+1=n-k-1$. By Facts 4 and 5, $G'\cong G_2$ or $G_4$. Note that $q_1(G_4)<q_1(G_2)$; hence by Lemma 2.5, our lemma holds in this case.

{\bf Case 2.}\ $l(G)>l_G(v_1,v_2,\ldots,v_k)$. Without loss of generality, we assume that $l(G)=l_G(v_2,v_3, \ldots,v_k,v_{k+1})=s$. Hence, $s< n-k-3$. Let $\delta_1,\delta_2,\ldots,\delta_s$ be the $(k+1)$-cliques satisfying the property $P_G(v_2,v_3,\ldots,v_k,v_{k+1})$ and $w_1,w_2,\ldots,w_s$ be the vertices in these $(k+1)$-cliques different from $v_2,v_3,\ldots,v_k, v_{k+1}$.

If $l_G(v_1,v_2,\ldots,v_k)\neq 0$, then by a similar discussion as in the proof of Case 1, one can show this lemma is true. Hence, we consider that $l_G(v_1,v_2,\ldots,v_k)=0$. By the definitions of the $k$-tree and $S_1(G)$, it is routine to check that there exists a vertex $r_1$ in $S_1(G)$ such that
$v_{k+2}r_1\in E_G$. Thus we have $v_{k+1}r_1\not\in E_G$ and there exists only one vertex, say $v_i$, in $\{v_1,v_2,\ldots,v_k\}$ such that $v_ir_1\not\in E_G$. If $i=1$, then $N_G(r_1)=\{v_2,v_3,\ldots,v_k,v_{k+2}\}$. By a similar discussion as in the proof of Case 1, one can show this lemma holds. In what follows we consider that $i\neq1$. Then, $N_G(r_1)=\{v_1,v_2,\ldots,v_{i-1},v_{i+1},\ldots,v_k,v_{k+2}\}$.

Assume that $l_G(v_1,v_2,\ldots,v_{i-1},v_{i+1},\ldots,v_k,v_{k+2})=t$, then $1\leqslant    t\leqslant    s$. Let $\Delta_1,\Delta_2,\ldots,\Delta_t$
be the $(k+1)$-cliques satisfying the property $P_G(v_1,v_2,\ldots,v_{i-1},v_{i+1},\ldots,v_k,v_{k+2})$ and $r_1,r_2,\ldots,r_t$ be the vertices in these $(k+1)$-cliques different from $v_1,v_2,\ldots,v_{i-1},v_{i+1},\ldots,v_k,v_{k+2}$. We proceed through the following two possible subcases.

\textbf{Subcase 2.1.}\ $x_{v_1}\leqslant x_{v_i}$. Then let $G'=G-v_1r_1+v_ir_1$. By Lemma 1.1, we have $q_1(G) < q_1(G')$. Note that $\{v_2,v_3,\ldots,v_k,v_{k+1},v_{k+2}\}\subseteq N_{G'}(v_1)$; hence we get $l(G')=l(G),S_1(G')=S_1(G)$ and $N_{G'}(r_1)=\{v_2,v_3,\ldots,v_k,\linebreak v_{k+2}\}.$ By a similar discussion as in the proof of Case 1, one can find an $n$-vertex $k$-tree $G^*$ with $l(G^*)=l(G')+1$ such that $q_1(G') < q_1(G^*)$ and either (2.1) or (2.2) holds. Hence, we have $q_1(G) < q_1(G') < q_1(G^*)$ and $l(G^*)=l(G')+1=l(G)+1$, as desired.

\textbf{Subcase 2.2.}\ $x_{v_1} > x_{v_i}$. Obviously, $G\ [\{v_1,v_2,\ldots,v_{i-1},v_{i+1},\ldots,v_k,v_{k+1}\}\ ]$ is a $k$-clique. Assume that $l_G(v_1,v_2,\linebreak \ldots, v_{i-1},v_{i+1},\ldots,v_k,v_{k+1})=a$ where $0\leqslant    a\leqslant    s$. Let $\Delta'_1,\Delta'_2,\ldots,\Delta'_a$ be the $(k+1)$-cliques satisfying the property $P_G(v_1,v_2,\ldots, v_{i-1},v_{i+1},\ldots,v_k,v_{k+1})$ and $y_1,y_2,\ldots,y_a$ be the vertices in these $(k+1)$-cliques different from $v_1,v_2,\ldots,v_{i-1},v_{i+1},\ldots,v_k,v_{k+1}$. Put
$$
  G'':=G-\{v_iw_1,v_iw_2,\ldots,v_iw_{s-a}\}+\{v_1w_1,v_1w_2,\ldots,v_1w_{s-a}\}.
$$
Then $G''$ is in $\mathscr{T}_n^k$ and by Lemma 1.1 we have $q_1(G) < q_1(G'')$. Note that $\{v_1,\ldots,v_{i-1},v_{i+1},\ldots,v_k,v_{k+1},v_{k+2}\}\subseteq N_{G''}(v_i)$; hence we get $l(G'')=l(G),S_1(G'')=S_1(G)$. On the other hand, we have
$$
  N_{G''}(y_t)=N_{G''}(w_j)=\{v_1,v_2,\ldots,v_{i-1},v_{i+1},\ldots,v_k,v_{k+1}\}
$$
for all $0\leqslant    t\leqslant    a$ and $1\leqslant    j\leqslant    s-a$. By a similar discussion as in the proof of Case 1, there exists an $n$-vertex $k$-tree $G^*$
with $l(G^*)=l(G'')+1$ such that $q_1(G'') < q_1(G^*)$ and either (2.1) or (2.2) holds.
Hence, we have $q_1(G) < q_1(G'') < q_1(G^*)$ and $l(G^*)=l(G'')+1=l(G)+1$, as desired.

This completes the proof.
\end{proof}

\begin{thm}
Let $G$ be in $\mathscr{T}_n^k\setminus \{S_{k,n-k},G_1\}$ on $n\geqslant
k+1$ vertices. Then, $q_1(G)\leqslant  q_1(G_2)$ and the equality holds if
and only if $G\cong G_2,$ where $G_2$ is depicted in Fig. 2.
\end{thm}
\begin{proof}
Since $G\ncong S_{k,n-k},\ G_1$, we have $|S_1(G)|\leqslant n-k-1$ and $l(G)\leqslant n-k-3$. Note that $|S_1(G_1)|=|S_1(G_2)|= n-k-1$; hence we proceed by considering the following possible cases.

\textbf{Case 1.}\ $|S_1(G)|= n-k-1$. In this case, let $G$ be in $\mathscr{T}_n^k\setminus \{G_1\}$ with $|S_1(G)|=n-k-1$ such that $q_1(G)$ is as large as possible. By Lemma 2.3, we have $l(G)=n-k-3$. By Facts 4 and 5 we have $G\cong G_2$ or $G_4$. In view of Lemma 2.5, $q_1(G_4)< q_1(G_2)$, one can easily get $q_1(G)\leqslant  q_1(G_2)$ when $G\in \mathscr{T}_n^k\setminus \{G_1\}$ and the equality holds if
and only if $G\cong G_2$.

\textbf{Case 2.} $|S_1(G)|= n-k-2$. In this case, by Fact 2 we obtain that $G-S_1(G)$ is a $k$-tree on $k+2$ vertices, that is, $G-S_1(G)\cong S_{k,2}$. Denote the vertex set of $G-S_1(G)$ by $\{v_1,v_2,\ldots,v_{k+2}\}$ with $d_{G-S_1(G)}(v_{k+1})=d_{G-S_1(G)}(v_{k+2})=k$. Obviously, $v_{k+1}v_{k+2}\not\in E_G$ and there exist vertices $w_1$ and $r_1$ in $S_1(G)$ such that $v_{k+1}w_1,\linebreak v_{k+2}r_1\in E_G$ by the definitions of the $k$-tree and $S_1(G)$. Thus, $N_G(r_1)\not= N_G(w_1).$  Hence, $l(G)<|S_1(G)|=n-k-2$, i.e., $l(G)\leqslant n-k-3$.

Let $G$ be in $\mathscr{T}_n^k$ with $|S_1(G)|=n-k-2$ such that $q_1(G)$ is as large as possible. By Lemma 2.6, one has either $l(G)=n-k-3$ or $q_1(G)< q_1(G_2)$. However, if $l(G)=n-k-3$, then combining with $|S_1(G)|=n-k-2$ and by Facts 4 and 5, we have $G\cong G_3$ or $G_5$. By Lemma 2.5, $q_1(G_5)<q_1(G_3)< q_1(G_2).$ Hence, we obtain $q_1(G) < q_1(G_2)$ if $|S_1(G)|= n-k-2$.

\textbf{Case 3.} $|S_1(G)|\leqslant    n-k-3$. In this case, repeatedly using Lemma 3.1 yields a $k$-tree $G'$ on $n$ vertices such that $q_1(G)<q_1(G')$ with $|S_1(G')|=n-k-2$. However, in view of Case 2, we get $q_1(G')<q_1(G_2)$. So, $q_1(G)<q_1(G')<q_1(G_2)$ if $G$ is an $n$-vertex $k$-tree with $|S_1(G)|\leqslant    n-k-3$.

By Cases 1-3, Theorem 2.7 holds.
\end{proof}

\section*{\normalsize Acknowledgements}
We thank Dr. Asghar Bahmani for drawing our attention to the mistake.

\end{document}